
\input amstex.tex
\documentstyle{amsppt}

\def\DJ{\leavevmode\setbox0=\hbox{D}\kern0pt\rlap
 {\kern.04em\raise.188\ht0\hbox{-}}D}
\def\dj{\leavevmode
 \setbox0=\hbox{d}\kern0pt\rlap{\kern.215em\raise.46\ht0\hbox{-}}d}

\def\txt#1{{\textstyle{#1}}}
\baselineskip=13pt
\def\hf{{\textstyle{1\over2}}}
\def\a{\alpha}
\def\d{{\,\roman d}}
\def\e{\varepsilon}
\def\f{\varphi}
\def\G{\Gamma}

\def\s{\sigma}

\def\={\;=\;}

\def\zt{\zeta(\hf+it)}

\def\no{\noindent}  
\def\R{\Re{\roman e}\,} \def\I{\Im{\roman m}\,} 
\def\z{\zeta}

\def\hf{{\textstyle{1\over2}}}
\def\txt#1{{\textstyle{#1}}}
\def\f{\varphi}

\font\tenmsb=msbm10
\font\sevenmsb=msbm7
\font\fivemsb=msbm5
\newfam\msbfam
\textfont\msbfam=\tenmsb
\scriptfont\msbfam=\sevenmsb
\scriptscriptfont\msbfam=\fivemsb
\def\Bbb#1{{\fam\msbfam #1}}

\def \NN {\Bbb N}

\def \RR {\Bbb R}

\font\ff=cmr8
\def\txt#1{{\textstyle{#1}}}
\baselineskip=13pt

\font\teneufm=eufm10
\font\seveneufm=eufm7
\font\fiveeufm=eufm5
\newfam\eufmfam
\textfont\eufmfam=\teneufm
\scriptfont\eufmfam=\seveneufm
\scriptscriptfont\eufmfam=\fiveeufm
\def\mathfrak#1{{\fam\eufmfam\relax#1}}

\font\tenmsb=msbm10
\font\sevenmsb=msbm7
\font\fivemsb=msbm5
\newfam\msbfam
     \textfont\msbfam=\tenmsb
      \scriptfont\msbfam=\sevenmsb
      \scriptscriptfont\msbfam=\fivemsb
\def\Bbb#1{{\fam\msbfam #1}}

\def \NN {\Bbb N}

\def \RR {\Bbb R}

  \def\rightheadline{{\hfil{\ff
  Some identities for the Riemann zeta-function}\hfil\tenrm\folio}}

  \def\leftheadline{{\tenrm\folio\hfil{\ff
   A. Ivi\'c }\hfil}}
  \def\emptyheadline{\hfil}
  \headline{\ifnum\pageno=1 \emptyheadline\else
  \ifodd\pageno \rightheadline \else \leftheadline\fi\fi}

\font\ff=cmr8
\font\teneufm=eufm10
\font\seveneufm=eufm7
\font\fiveeufm=eufm5
\newfam\eufmfam
\textfont\eufmfam=\teneufm
\scriptfont\eufmfam=\seveneufm
\scriptscriptfont\eufmfam=\fiveeufm
\def\mathfrak#1{{\fam\eufmfam\relax#1}}

\font\tenmsb=msbm10
\font\sevenmsb=msbm7
\font\fivemsb=msbm5
\newfam\msbfam
\textfont\msbfam=\tenmsb
\scriptfont\msbfam=\sevenmsb
\scriptscriptfont\msbfam=\fivemsb
\def\Bbb#1{{\fam\msbfam #1}}

\def \NN {\Bbb N}

\def \RR {\Bbb R}

\def\a{\alpha} 
 \def\e{\varepsilon}
\def\no{\noindent} \def\d{\,{\roman d}}
\topmatter
\title
Some  identities for the Riemann zeta-function
\endtitle
\author
Aleksandar Ivi\'c
\endauthor
\address
Katedra Matematike RGF-a Universiteta u Beogradu, \DJ u\v sina 7,
11000 Beograd, Serbia (Yugoslavia)
\endaddress
\subjclass 11 M 06
\endsubjclass
\email {\tt eivica\@ubbg.etf.bg.ac.yu, aivic\@rgf.bg.ac.yu}
\endemail
\dedicatory
To appear in ``Publikac. Elektrotehn. Fak. Ser. Mat." 2004
\enddedicatory
\abstract
Several identities for the Riemann zeta-function
$\zeta(s)$ are proved. For example, if $s = \s + it$ and $\s > 0$, then
$$
\int_{-\infty}^\infty \left|{(1-2^{1-s})\z(s)\over s}\right|^2\d t
= {\pi\over\s}(1 - 2^{1-2\s})\z(2\s).
$$
\endabstract
\endtopmatter

Let as usual $\z(s) = \sum_{n=1}^\infty n^{-s}\;(\R s > 1)$ denote
the Riemann zeta-function. The motivation for this note is the
quest to evaluate explicitly integrals of $|\zt|^{2k},\,k\in\NN$,
weighted by suitable functions. In particular, the problem is
to evaluate  in closed form
$$
\int_0^\infty (3-\sqrt{8}\cos(t\log2))^k|\zt|^{2k}
{\d t\over{(\txt{1\over4}}+t^2)^k} \qquad(k\in\NN).
$$
When $k = 1,2$ this may be done, thanks to the identities which will
be established below. The first identity in question is given by

\bigskip
THEOREM 1. {\it Let $s = \s + it$. Then for $\s > 0$ we have}
$$
\int_{-\infty}^\infty \left|{(1-2^{1-s})\z(s)\over s}\right|^2\d t
= {\pi\over\s}(1 - 2^{1-2\s})\z(2\s).\leqno(1)
$$

\bigskip\no
Since $\lim_{s\to1}(s-1)\z(s) =1 $, then setting in (1) $\s = \hf$ 
we obtain the following

\bigskip
{\bf Corollary 1}.
$$
\int_0^\infty (3-\sqrt{8}\cos(t\log2))|\zt|^2{\d t\over{\txt{1\over4}}+t^2}
\;=\; \pi\log2.\leqno(2)
$$

\bigskip
Another identity, which relates directly the square  of $\z(s)$ to
a Mellin-type integral, is contained in

\bigskip
THEOREM 2. {\it Let $\chi_{\Cal A}(x)$ denote the characteristic function
of the set $\Cal A$, and let}
$$
\f(x) := \sum_{m=1}^\infty \sum_{n=1}^\infty \int_1^x
\chi_{[2m-1,2m)}({x\over u})\chi_{[2n-1,2n)}(u)\,{\d u\over u}
\qquad(x \geqslant 1).
\leqno(3)
$$
{\it Then for $\s > 0$ we have}
$$
s^2\int_1^\infty \f(x)x^{-s-1}\d x = (1 - 2^{1-s})^2\z^2(s).\leqno(4)
$$

\bigskip\no
From (4) we obtain the following

\medskip
{\bf Corollary 2}.
$$
\int_0^\infty (3-\sqrt{8}\cos(t\log2))^2|\zt|^4
{\d t\over({\txt{1\over4}}+t^2)^2}
\;=\; \pi\int_1^\infty \f^2(x){\d x\over x^2}.\leqno(5)
$$
The integral on the right-hand side of (5) is elementary, but
nevertheless its evaluation in closed form is complicated.

\bigskip
{\bf Proof of Theorem 1}. We start from (see e.g., [1, Chapter 1])
the identity
$$
(1 - 2^{1-s})\z(s) \;=\; \sum_{n=1}^\infty (-1)^{n-1}n^{-s}\qquad(\s > 0)
\leqno(6)
$$
and
$$
\int_{-\infty}^\infty{\cos(\a x)\over\s^2 + x^2}\,\d x \;=\; 
{\pi\over\s}e^{-|\a|\s}\qquad(\a\in\RR,\;\s>0),\leqno(7)
$$
which follows by the residue theorem on integrating $e^{i\a z}/(\s^2+z^2)$
over the contour consisting of $[-R,\,R]$ and semicircle $|z| = R, \I z > 0$
and letting $R\to\infty$. By using (6) and (7) it is seen that the
left-hand side of (1) becomes
$$
\eqalign{&
\sum_{m=1}^\infty\sum_{n=1}^\infty(-1)^{m+n}(mn)^{-\s}\int_{-\infty}^\infty
\left({m\over n}\right)^{it}{\d t\over\s^2+t^2}\cr&
= {\pi\over\s}\z(2\s) + 
2\sum_{m=1}^\infty\sum_{n<m}(-1)^{m+n}(mn)^{-\s}\int_{-\infty}^\infty
{\cos(t\log{m\over n})\over\s^2+t^2}\,\d t\cr&
= {\pi\over\s}\left(\z(2\s) + 2\sum_{m=1}^\infty(-1)^mm^{-\s}\sum_{n=1}^{m-1}
(-1)^nn^{-\s}\cdot e^{-\s\log{m\over n}}\right)\cr&
= {\pi\over\s}\left(\z(2\s) + 2\sum_{m=1}^\infty(-1)^mm^{-2\s}
\sum_{n=1}^{m-1}(-1)^n\right)\cr&
= {\pi\over\s}\left(\z(2\s) + 2\sum_{k=1}^\infty(-1)^{2k}(2k)^{-2\s}(-1)
\right) = {\pi\over\s}(1-2^{1-2\s})\z(2\s).\cr}
$$
This holds initially for $\s>1$, but by analytic continuation it holds for
$\s>0$ as well.
\bigskip\no
We shall provide now a second proof of Theorem 1. As in the formulation of
Theorem 2, let $\chi_{\Cal A}(x)$ denote the characteristic
function of the set ${\Cal A}$, and let the interval $[a, b)$ denote the
set of numbers $\{x : a \leqslant x < b\}$. Then, for $\s > 0$, we have
$$\eqalign{&
\int_1^\infty x^{-s-1}\sum_{n=1}^\infty \chi_{[2n-1,2n)}(x)\d x 
= \sum_{n=1}^\infty \int_{2n-1}^{2n}x^{-s-1}\d x\cr&
= {1\over s}\sum_{n=1}^\infty \left((2n-1)^{-s}-(2n)^{-s}\right)
= {(1-2^{1-s})\z(s)\over s}\cr}\leqno(8)
$$
in view of (6).
Now we invoke Parseval's identity for Mellin transforms (see e.g., [1]
and [3]). We need this identity for the modified
 Mellin transforms, defined by
$$
F^*(s) \equiv m[f(x)] := \int_1^\infty f(x)x^{-s-1}\d x.
$$
The properties of this transform were developed by the author in [2]. 
In particular, we need Lemma 3 of [2] which says that
$$
\int_1^\infty f(x)g(x)x^{1-2\s}\d x =  {1\over2\pi i}\int_{\R s = \s}
F^*(s)\overline{G^*(s)}\d s \leqno(9)
$$
if $F^*(s) = m[f(x)],\,G^*(s) = m[g(x)]$, and $f(x), g(x)$ are
real-valued, continuous functions for $x > 1$, such that
$$
x^{{1\over2}-\s}f(x) \in L^2(1,\infty),\quad x^{{1\over2}-\s}g(x)
\in L^2(1,\infty).
$$
From (8) and (9) we obtain, for $\s>0$,
$$
\int_1^\infty{1\over x^2}\left(\sum_{n=1}^\infty
\chi_{[2n-1,2n)}(x)\right)^2x^{1-2\s}\d x
= {1\over2\pi i}\int_{\R s = \s}\left|{(1-2^{1-s})\z(s)\over s}\right|^2\d s.
$$
But as $\chi_{\Cal A}^2(x) = \chi_{\Cal A}(x)$, it is easily found that
the left-hand side of the above identity equals
$$\eqalign{&
\sum_{m=1}^\infty\sum_{n=1}^\infty \int_1^\infty
\chi_{[2m-1,2m)}(x)\chi_{[2n-1,2n)}(x)x^{-1-2\s}\d x\cr&
= \sum_{n=1}^\infty  \int_{2n-1}^{2n}x^{-1-2\s}\d x = {(1-2^{1-2\s})\z(2\s)
\over2\s}\cr}
$$
in view of (6), and (1) follows. 

\bigskip
For the Proof of Theorem 2 we need the following

\medskip
LEMMA. {\it Let $0 < a < b$. If $f(x)$ is integrable on $[a,b]$, then}
$$\eqalign{
&{\left(\int_a^b f(x)x^{-s}\d x\right)}^2\cr&
= \int_{a^2}^{ab}x^{-s}\int_a^{x/a}f(u)f({x\over u})
{\d u\over u}\d x
+ \int_{ab}^{b^2}x^{-s}\int_{x/b}^{b}f(u)f({x\over u})
{\d u\over u}\d x.
\cr}\leqno(10)
$$
{\it The identity} (10) {\it remains valid if $b = \infty$, provided
the integrals in question converge, in which case the second integral
on the right-hand side is to be omitted.}

\medskip
{\bf Proof of the Lemma.} We write the left-hand side of (10) as
the double integral
$$
\int_a^b\int_a^b(xy)^{-s}f(x)f(y)\d x\d y
$$
and make the change of variables $x = X/Y, y = Y$.
The Jacobian of this transformation
equals $1/Y$, hence the left-hand side of (10) becomes
$$\eqalign{&
\int_{a^2}^{b^2}X^{-s}\left(\int_{\max(a,X/b)}
^{\min(X/a,b)} f(Y)f({X\over Y})\,
{\d Y\over Y}\right)\d X\cr&
= \int_{a^2}^{ab}X^{-s} \int_a^{X/a}f(Y)f({X\over Y})\,
{\d Y\over Y}\d X\cr&
+ \int_{a^b}^{b^2}X^{-s} \int_{X/b}^{b}f(Y)f({X\over Y})\,
{\d Y\over Y}\d X,\cr&}
$$
as asserted.

\bigskip
{\bf Proof of Theorem 2.}  We use (8) and the Lemma to obtain that (4)
certainly holds with $\f(x)$ given by (3), since trivially $\f(x) \ll x$.
To see that it holds for $\s>0$,
we note that 
$$
\int_1^x g(u)g({x\over u}){\d u\over u} =
\int_1^{\sqrt x} + \int_{\sqrt x}^x = 2\int_{\sqrt x}^x 
g(u)g({x\over u}){\d u\over u},\leqno(11)
$$
and use (11) with
$$
g(x) = \sum_{n=1}^\infty \chi_{[2n-1,2n)}(x).
$$
Note then that the integrand
in $\f(x)$ equals $1/u$ for $2m - 1 \leqslant u \leqslant 
2m,\, 2n-1\leqslant u \leqslant 2n$,
and otherwise it is zero. This gives the condition
$$
4mn - 2m - 2n + 1 \leqslant x < 4mn, \hf\sqrt{x} \leqslant n \le \hf(x+1),\;
1 \leqslant m \leqslant \hf(\sqrt{x}+1).
$$
We also have
$$
\int_{\sqrt{x}}^x \chi_{[2m-1,2m)}({x\over u})\chi_{[2n-1,2n)}(u){\d u
\over u} \leqslant \int_{2n-1}^{2n}{\d u\over u} \leqslant {1\over 2n-1}.
$$
Therefore
$$\eqalign{
\f(x) & \ll \sum_{m\leqslant\sqrt{x}}\;\sum_{x/(4m)<n\leqslant(x-1+2m)
/(4m-2)}{1\over n}
\cr&
\ll \sum_{m\leqslant\sqrt{x}}{m\over x}
\left(1 + {x\over m^2}\right) \ll \log x.
\cr}\leqno(12)
$$
This bound shows that the integral in (4) is absolutely convergent
for $\s > 0$. Thus by the principle of analytic continuation this
 completes the proof of  Theorem 2.

\medskip Corollary 2 follows then from (4) and (9) on setting $\s = \hf$.

\medskip
It is interesting to note that the bound in (12) is actually of the
correct order of magnitude. Namely we have

\bigskip
THEOREM 3. {\it For any given $\e > 0$ we have}
$$
\f(x) \,=\, {\txt{1\over4}}\log x + {\txt{1\over2}}\log\left(
{\txt{\pi\over2}}\right) +  O_\e\left(x^{\e-{1\over4}}\right).
\leqno(13)
$$

\bigskip {\bf Proof of Theorem 3.} By (8) and the inversion formula for
the Mellin transform $m[f(x]$ (see [2, Lemma 1]) we have, for any
$c > 0$,
$$
\f(x) = {1\over2\pi i}\int_{\R s = c}{(1-2^{1-s})^2\z^2(s)x^s\over s^2}\d s.
\leqno(14)
$$
We shift the line of integration in (14) to $c = \e - 1/4$ with 
$0 < \e < 1/8$, which clearly may be assumed. Since $\z(0) = -\hf$
and $\z'(0) = -\hf\log(2\pi)$, the
residue at the double pole $s=0$ is found to be
$$
{\txt{1\over4}}\log x + A,\quad A = -\z'(0) - \log 2
= \hf\log\left({\txt{\pi\over2}}\right).\leqno(15)
$$
We use the functional equation (see e.g., [1, Chapter 1]) for $\z(s)$,
namely
$$
\z(s) = \chi(s)\z(1 - s),\quad \chi(s) = 2^s\pi^{s-1}\sin(\hf\pi s)
\G(1-s)
$$
with
$$
\chi(s) = {\left({2\pi\over t}\right)}^{\s+it-{1\over2}}{\roman e}
^{i(t+{1\over4}\pi)}\cdot\left(1 + O\left({1\over t}\right)\right)
\quad(t \geqslant 2).
$$
Let $s = \e-{1\over4}+it$. Then by absolute convergence we have
$$
\eqalign{&
\int_T^{2T}
{(1-2^{1-s})^2\z^2(s)x^s\over s^2}\d t\cr&
= i\sum_{n=1}^\infty d(n)n^{\e-5/4}\int_T^{2T}{(1-2^{1-s})^2\over s^2}
x^{\e-{1\over4}+it}\left({t\over2\pi}\right)^{{3\over2}-2\e}{\roman e}
^{iF(t,n)}\d t + O(T^{-{1\over2}-2\e}),\cr}
$$
where $d(n)$ is the number of divisors of $n$ and
$$
F(t,n) := 2t + t\log n - 2t\log(t/2\pi),\quad
{\d^2\over \d t^2}\,(t\log x + F(t,n)) = -{2\over t}\,.
$$
Hence by the second derivative test (see [1, Lemma 2.2]) the 
above series is 
$$
\ll \sum_{n=1}^\infty d(n)n^{\e-5/4}T^{-2\e} 
= \z^2({\txt{5\over4}}-2\e)T^{-2\e}
\ll T^{-2\e}.
$$
This shows that
$$
\int_{\R s = \e-1/4}{(1-2^{1-s})^2\z^2(s)x^s\over s^2}\d s \ll x^{\e-1/4},
$$
hence (13) follows from (14), (15) and the residue theorem.

\medskip
In concluding, note that if we write
$$
\f(x) = {\txt{1\over4}}\log x + A + \f_1(x),
$$
where $A$ is given by (15) then, for $\R s = \s  > 0$, (4) yields
$$
s^2\left({A\over s} + {1\over 4s^2} + \int_1^\infty \f_1(x)x^{-s-1}
\d x\right) = (1 - 2^{1-s})^2\z^2(s),
$$
and the above integral converges absolutely, for $\s > -1/4$, 
in view of (13). Thus by
analytic continuation it follows that, for $\s > -1/4$,
$$
As + {\txt{1\over4}} + s^2\int_1^\infty \f_1(x)x^{-s-1}\d x 
= (1 - 2^{1-s})^2\z^2(s).
$$


\bigskip\bigskip
\Refs
\bigskip

\item{[1]} A. Ivi\'c, The Riemann zeta-function, John Wiley \&
Sons, New York, 1985.

\item{[2]} A. Ivi\'c,  On some conjectures and results for the
Riemann zeta-function and Hecke series, Acta Arith. {\bf109}(2001),
115-145.

\item{[3]} E.C. Titchmarsh, Introduction to the Theory of Fourier
Integrals, Oxford University Press, Oxford, 1948.

\bigskip

Aleksandar Ivi\'c

Katedra Matematike RGF-a

Universitet u Beogradu

\DJ u\v sina 7, 11000 Beograd, Serbia

{\tt aivic\@rgf.bg.ac.yu}

\endRefs


\bye